\documentclass[12pt]{article}

\newcommand{\be}{\begin{equation}}
\newcommand{\ee}{\end{equation}}
\newtheorem{mrem}{Remark}
\newtheorem{mthm}{Theorem}

\newtheorem{pro}{Proposition}[section]

\begin{document}

\title{ A Liouville Theorem for the Fractional Laplacian}
\author {Ran Zhuo \,\, Wenxiong Chen \,\, Xuewei Cui \,\, Zixia Yuan}

\date{}
\maketitle

\begin{abstract}

In this paper, we consider the following fractional Laplace equation
\begin{eqnarray}\left\{\begin{array}{ll}
(- \Delta)^{\alpha/2}u(x)=0,& \mbox {in\,\,\,} R^n,\\
u(x)\geq 0,& \mbox { in\,\,\,} R^n,
\end{array} \right.\label{a1} \end{eqnarray}
where $n \geq 2$ and $\alpha$ is any real number between $0$ and $2$. We prove that the only solution
for (\ref{a1}) is constant. Or equivalently,

{\em Every $\alpha$-harmonic function bounded either above or below in all of $R^n$ must be constant.}

This extends the classical Liouville Theorem from
Laplacian to the fractional Laplacian.

As an immediate application, we use it to obtain an equivalence between a semi-linear pseudo-differential
\begin{equation}
(-\Delta)^{\alpha/2} u = u^p (x) ,  \;\; x \in R^n
\label{b1}
\end{equation}
and the corresponding integral equation
$$ u(x) = \int_{R^n} \frac{1}{|x-y|^{n-\alpha}} d y , \;\; x \in R^n. $$ Combining this with the existing results on
the integral equation, one can obtained much more general results on the qualitative properties of the solutions for (\ref{b1}).

\end{abstract}

\maketitle

\section{Introduction}

The well-known Liouville's Theorem states that

{\em Any harmonic function bounded either above or below in all of $R^n$ is constant.}

One of its important applications is the proof of the Fundamental Theorem of Algebra. It is
also a key ingredient in deriving a priori estimates for solutions in PDE analysis.

The main purpose of this article is to extend this classical theorem to the fractional Laplacian.

Essentially different from the Laplacian, the fractional Laplacian in $R^n$ is a nonlocal operator, taking the form
\begin{equation}
(-\Delta)^{\alpha/2} u(x) = C_{n,\alpha} \, PV \int_{R^n} \frac{u(x)-u(z)}{|x-z|^{n+\alpha}} dz
\label{Ad7}
\end{equation}
 where $\alpha$ is any real number between $0$ and $2$ and PV stands for the Cauchy principal value.
This operator is well defined in $\cal{S}$, the Schwartz space of rapidly decreasing $C^{\infty}$
functions in $R^n$. In this space, it can also be defined equivalently in terms of the Fourier transform
$$ \widehat{(-\Delta)^{\alpha/2} u} (\xi) = |\xi|^{\alpha} \hat{u} (\xi) $$
where $\hat{u}$ is the Fourier transform of $u$. One can extend this operator to a wider space of distributions as the following.

Let
$$L_{\alpha}=\{u: R^n\rightarrow R \mid \int_{R^n}\frac{|u(x)|}{(1+|x|^{n+\alpha})}<\infty\}.$$

For $u\in L_{\alpha}$, we define $(-\Delta)^{\alpha/2} u$ as a distribution:
$$< (-\Delta)^{\alpha/2}u(x), \phi> \, = \, < u, (-\Delta)^{\alpha/2} \phi > ,  \;\;\; \forall \, \phi \in C_0^{\infty} (R^n) . $$
This defines the operator in a weak sense. In this paper, we consider the class of functions where the fractional Laplacian is defined in a little bit more stronger sense, that is

$u \in L_{\alpha}$,
such that the right hand side of (\ref{Ad7}) is well defined for every $x \in R^n$.

One can verify that, all the above definitions coincides when $u$ is in $\cal{S}$.

We say that $u$ is an $\alpha$-harmonic function if $u \in L_{\alpha}$,
such that the right hand side of (\ref{Ad7}) is well defined for every $x \in R^n$ and equals zero. In this sense, we have
\begin{mthm}
Every $\alpha$-harmonic function bounded either above or below in all of $R^n$ for $n \geq 2$ must be constant.
\label{mthm1}
\end{mthm}
This is the main result of the paper. To prove it, we study
\begin{eqnarray}\left\{\begin{array}{ll}
(- \Delta)^{\alpha/2}u(x)=0,& \mbox {in\,\,\,} R^n,\\
u(x)\geq 0,& \mbox { in\,\,\,} R^n.
\end{array} \right.\label{a}
\end{eqnarray}
We say that $u \geq 0$ is a strong solution of (\ref{a}), if $u \in L_{\alpha}$,
such that the right hand side of (\ref{Ad7}) is well defined for every $x \in R^n$ and equals zero.

Apparently, Theorem \ref{mthm1} is equivalent to the following
\begin{mthm}
Assume that $n \geq 2$. Let $u$ be a strong solution of (\ref{a}), then $u \equiv C$.
\label{mthm2}
\end{mthm}

As an immediate application of the Liouville theorem for $\alpha$ harmonic functions,
we prove an equivalence between a psedodifferential equation and an integral equation.

\begin{mthm}
Assume that $n \geq 2$ and $u \in L_{\alpha}$ is a nonnegative strong solution of
\begin{equation}
(- \Delta)^{\alpha/2}u(x)=u^p(x), \;\; x \in R^n,
\label{3-1}
\end{equation}
then $u$ also satisfies
$$u(x)=\int_{R^n}\frac{c_n}{|x-y|^{n-\alpha}}u^p(y)dy,$$
and vice versa.
\label{mthm3}
\end{mthm}

\begin{mrem}

i) Actually, the right hand side of equation (\ref{3-1}) can be a much more general function $f(x, u)$,
such that, for any constant $c>0$,
\begin{equation}
\int_{R^n}\frac{1}{|x-y|^{n-\alpha}} f(y, c) dy = \infty.
\label{z1}
\end{equation}

ii) The idea of proof can be extended to establish the equivalence between a general system of
m equations in $R^n$
$$\left\{\begin{array}{ll}
(-\Delta)^{\alpha/2} u_i (x)= f_i(x, u_1(x), \cdots u_m(x)), & i=1, \cdots, m, \\
u_i \geq 0 , & i=1, \cdots, m,
\end{array}
\right.
$$ and the corresponding integral system
$$ \left\{\begin{array}{ll}
u_i(x) = \int_{R^n} \frac{c_n}{|x-y|^{n-\alpha}} f_i(y, u_1(y), \cdots, u_m(y)), & i=1, \cdots, m, \\
u_i(x) \geq 0, & i=1, \cdots, m.
\end{array}
\right.
$$
\end{mrem}

Combining Theorem \ref{mthm3} with the qualitative properties established for the integral equations
in \cite{CLO} and \cite{CLO1}, one obtain immediately that
\begin{mthm}
Assume that $n \geq 2$ and $u$ is a nonnegative strong solution of (\ref{3-1}) for $0<\alpha < 2$. Then

i) In  the critical case when $p = \frac{n+\alpha}{n-\alpha}$, it must assume the form
$$ u(x) = c (\frac{t}{t^2 + |x -
x_o|^2})^{(n-\alpha)/2}$$
for some  $t>0$, $x_o\in{\mathbf{R}}^n$.

ii) In the subcritical case when $1<p<\frac{n+\alpha}{n-\alpha}$, we must have $u\equiv 0.$
\label{mthm4}
\end{mthm}

\begin{mrem}

i) In \cite{CLO} and \cite{CLO1}, in order the results in Theorem \ref{mthm4} to hold, one requires
$u$ to be in $H^{\alpha/2}(R^n)$. Here we only requires $u \in L_{\alpha}$, a much weaker restriction.

ii) In \cite{BCPS}, by using the extension method to obtain the same results as in Theorem \ref{mthm4}, the authors require that $1\leq \alpha < 2$ and $u$ be bounded. Obviously, our condition here is much weaker.

\end{mrem}

In Section 2, we prove the Liouville Theorem \ref{mthm2} and hence Theorem \ref{mthm1}. In Section 3,
we establish the equivalence and hence prove Theorem \ref{mthm3} and \ref{mthm4}.

\section{The proof of the Liouville Theorem}

In this section, we prove Theorem \ref{mthm2}.

{\bf{Proof.}} First, we define

\begin{eqnarray}u_k(x)=\left\{\begin{array}{ll}
u(x),& |x|\leq k,\\
\int_{B_k} P_k (y,x)u(y)dy,& |x|>k,
\end{array} \right.\label{1} \end{eqnarray}
where $P_k(y,x)$ is a Poisson kernel in the
exterior of the ball $B_k$ with radius $k$ and centered at the
origin:
\begin{equation} P_k (y,x)=\Gamma(\frac{n}{2})\pi^{-\frac{n}{2}-1}\sin(\frac{\pi
\alpha}{2})\frac{(|x|^2-k^2)^{\alpha/2}}
{(k^2-|y|^2)^{\alpha/2}}\frac{1}{|x-y|^n}, \;\; \;
|y|<k,\,\,|x|>k.\label{2}
\end{equation}

Obviously, for each $x$,
$$\lim_{k\rightarrow \infty}u_k(x)=u(x).$$

One can also verify (see \cite{L}) that
 \begin{equation}
 (-\Delta)^{\alpha/2}u_k (x) = 0 , \; \mbox{ for } |x|>k , \;\; \mbox{ and } u_k(x) \leq u(x) , \; \forall\, x \in R^n .
 \label{add1}
 \end{equation}
Moreover, by Taylor expansion, it is easy to derive that
\begin{equation}
u_k(x)=\frac{c_1}{|x|^{n-\alpha}}+O(\frac{1}{|x|^{n-\alpha+1}}).
\label{3}
\end{equation}

In order to prove that $u$ is constant, it suffice to show that for any $\psi\in C^\infty_0(R^n)$, satisfying the condition
\begin{equation}
\int_{R^n}\psi(x)dx=0 ,  \label{4}
\end{equation}
we have
$$\int_{R^n}u(x)\psi(x)dx=0.$$
Actually, we only need to prove
$$\lim_{k\rightarrow \infty}\int_{R^n}u_k(x)\psi(x)dx=0.$$

We divided the proof into two steps.

{\bf{Step 1.}} Let
$$\varphi(x)=\int_{R^n}\frac{\psi(y)}{|x-y|^{n-\alpha}}.$$

Combining Taylor expansion with (\ref{4}), we deduce that
\begin{equation}
\varphi(x) = O (\frac{1}{|x|^{n-\alpha+1}}) \,,\,\,\,as\,
|x|\rightarrow \infty.\label{7}
\end{equation}

It follows that  $\varphi(x)\in
L^2(R^n)$ for $n\geq 2$ , and
\begin{equation}
(-\Delta)^{\alpha/2}\varphi(x)=\psi(x)\,,\,\,\,x\in R^n.
\label{5}
\end{equation}

In this step, we will show that
\begin{equation}
\int_{R^n}u_k(x) \, (-\Delta)^{\alpha/2} \varphi(x)dx=
\int_{R^n}(-\Delta)^{\alpha/2}u_k(x) \, \varphi(x)dx.
\label{6}
\end{equation}

Let $$v_k(x)=u_k(x)-\frac{c_1}{|x|^{n-\alpha}},$$  then $v_k(x)\in L^2(R^n)$ due to (\ref{3}).

Applying the Parseval's formula to one part of the left hand side of (\ref{6}), we
derive that
\begin{eqnarray}
\int_{R^n}u_k(x)(-\Delta)^{\alpha/2}\varphi(x)dx
&=&\int_{R^n}(v_k(x)+\frac{c_1}{|x|^{n-\alpha}})(-\Delta)^{\alpha/2}\varphi(x)dx \nonumber\\
&=&\int_{R^n}v_k(x)(-\Delta)^{\alpha/2}\varphi(x)dx+\int_{R^n}\frac{c_1}{|x|^{n-\alpha}}(-\Delta)^{\alpha/2}\varphi(x)dx\nonumber\\
&=&\int_{R^n}\widehat{v_k}(\xi)|\xi|^{\alpha}  \overline{\widehat{\varphi}(\xi)} +c\varphi(0).
\label{9}
\end{eqnarray}
Here we have used a result in \cite{L} that, in the sense of distributions, the Fourier transform of $\frac{c_1}{|x|^{n-\alpha}}$ is a constant multiple of $|\xi|^{-\alpha}$, and by the definition of
Fourier transform on distributions (see \cite{L}), we have
$$
\int_{R^n}\frac{c_1}{|x|^{n-\alpha}}(-\Delta)^{\alpha/2}\varphi(x)dx = \int_{R^n} |\xi|^{-\alpha}
\overline{|\xi|^{\alpha} \widehat{\varphi}(\xi)} d \xi = c \varphi (0).$$
Also note that we are not able to apply the Parseval's formula directly to $\int_{R^n}u_k(x) \, (-\Delta)^{\alpha/2}\varphi(x)dx$ because $u_k$
may not be in $L^2(R^n)$.

For the right hand side of (\ref{6}), we have
\begin{eqnarray}
\int_{R^n}(-\Delta)^{\alpha/2}u_k(x) \; \varphi(x)dx
&=&\int_{R^n}(-\Delta)^{\alpha/2}(v_k(x)+\frac{c_1}{|x|^{n-\alpha}})\varphi(x)dx \nonumber\\
&=&\int_{R^n}(-\Delta)^{\alpha/2}v_k(x) \; \varphi(x)dx+\int_{R^n}(-\Delta)^{\alpha/2}(\frac{c_1}{|x|^{n-\alpha}}) \; \varphi(x)dx\nonumber\\
&=&\int_{R^n}|\xi|^{\alpha}\widehat{v_k}(\xi) \overline{\widehat{\varphi}(\xi)} d\xi+c\varphi(0).
\label{10}
\end{eqnarray}
Here we have used a well-known fact that $\frac{1}{|x|^{n-\alpha}}$ is a constant multiple of the fundamental solution of $(-\Delta)^{\alpha/2}$.


Now from (\ref{9}) and (\ref{10}), we arrive at (\ref{6}).
\medskip

{\bf{Step 2.}} We prove
\begin{equation}
\int_{R^n}(-\Delta)^{\alpha/2}u_k(x) \; \varphi(x)dx\rightarrow
0\,,\,\,\,as\, k\rightarrow \infty.\label{11}
\end{equation}

By elementary calculation, we separate the integral in (\ref{11})
into two parts,
\begin{eqnarray*}
\int_{R^n}(-\Delta)^{\alpha/2}u_k(x)\varphi(x)dx
&=&c\int_{R^n}\int_{R^n}\frac{u_k(x)-u_k(y)}{|x-y|^{n+\alpha}}dy\varphi(x)dx\\
&=&c\int_{B_r(0)}\int_{R^n}\frac{u_k(x)-u_k(y)}{|x-y|^{n+\alpha}}dy\varphi(x)dx\\
&+&
c\int_{R^n\setminus B_r(0)}\int_{R^n}\frac{u_k(x)-u_k(y)}{|x-y|^{n+\alpha}}dy\varphi(x)dx\\
&=&I_1+I_2,
\end{eqnarray*}
where $r<k$.

First, we consider $I_1$.

From the first equation of (\ref{a}), we have
\begin{eqnarray*}
0
&=&c\int_{R^n}\frac{u(x)-u(y)}{|x-y|^{n+\alpha}}dy \\
&=&c\int_{B_k(0)}\frac{u(x)-u(y)}{|x-y|^{n+\alpha}}dy +c\int_{R^n\setminus B_k(0)}\frac{u(x)-u(y)}{|x-y|^{n+\alpha}}dy.
\end{eqnarray*}
It follows that
\begin{eqnarray}
I_1&=&c\int_{B_r(0)}\int_{R^n}\frac{u_k(x)-u_k(y)}{|x-y|^{n+\alpha}}dy\varphi(x)dx \nonumber\\
&=&c\int_{B_r(0)}\int_{B_k(0)}\frac{u(x)-u_k(y)}{|x-y|^{n+\alpha}}dy\varphi(x)dx \nonumber\\
&+&
c\int_{B_r(0)}\int_{R^n\setminus B_k(0)}\frac{u(x)-u_k(y)}{|x-y|^{n+\alpha}}dy\varphi(x)dx \nonumber\\
&=& c\int_{B_r(0)}\int_{R^n\setminus
B_k(0)}\frac{u(y)-u_k(y)}{|x-y|^{n+\alpha}}dy\varphi(x)dx.
\label{14}
\end{eqnarray}

Applying (\ref{14}), we obtain, for each fixed $r$,
\begin{eqnarray}
|I_1|&\leq & c\int_{B_r(0)}\int_{R^n\setminus
B_k(0)}\frac{|u(y)|}{|x-y|^{n+\alpha}}dy|\varphi(x)|dx \nonumber\\
&\leq & c\int_{B_r(0)}\int_{R^n\setminus
B_k(0)}\frac{|u(y)|}{(1+|y|)^{n+\alpha}}dy|\varphi(x)|dx\rightarrow
0\,,\,\,\,as\,k\rightarrow \infty. \label{15}
\end{eqnarray}
Here we have used the facts that $u \in L_{\alpha}$ and $0 \leq u_k \leq u $.

Next we will show that
\begin{eqnarray}
I_2&=&c\int_{R^n\setminus
B_r(0)}(- \Delta)^{\alpha/2} u_k \, \varphi(x)dx
\rightarrow 0\,,\,\,\,as\,r\rightarrow \infty,  \label{16}
\end{eqnarray}
uniformly in $k$.

Let $$f_k(x)=(-\Delta)^{\alpha/2}u_k(x).$$

We first show that
\begin{equation}
\int_{R^n}\frac{f_k(y)}{|x-y|^{n-\alpha}}dy = u_k(x).  \label{17}
\end{equation}

To this end, denote $$g_k(x)=\int_{R^n}\frac{f_k(y)}{|x-y|^{n-\alpha}}dy.$$
For any $\psi(x)\in C^\infty_0(R^n)$, we show that
\begin{equation}
\int_{R^n}g_k(x)\psi(x)dx=\int_{R^n}u_k(x)\psi(x)dx. \label{18}
\end{equation}

Actually,
\begin{eqnarray*}
\int_{R^n}g_k(x)\psi(x)dx &=& \int_{R^n} \int_{R^n} \frac{(-\Delta)^{\alpha/2} u_k(y)}{|x-y|^{n-\alpha}} d y \, \psi (x) d x \\
&=& \int_{R^n} \int_{R^n} \frac{(-\Delta)^{\alpha/2} ( v_k(y)+ c_1/|y|^{n-\alpha})}{|x-y|^{n-\alpha}} d y \, \psi (x) d x\\
&=& \int_{R^n} \int_{R^n} \frac{(-\Delta)^{\alpha/2} v_k(y)}{|x-y|^{n-\alpha}} d y \, \psi (x) d x +
\int_{R^n} \frac{c_1}{|x|^{n-\alpha}} \psi (x) d x \\
&=& I_1 + I_2.
\end{eqnarray*}

Since both $(-\Delta)^{\alpha/2} v_k(x)$ and $\psi (x)$ are compactly supported, one can exchange the order of integration to derive
\begin{eqnarray*}
I_1 &=& \int_{R^n} (-\Delta)^{\alpha/2} v_k(y) \int_{R^n} \frac{\psi(x)}{|x-y|^{n-\alpha}} d x \, d y\\
&=& \int_{R^n} |\xi|^{\alpha} \widehat{v_k} (\xi) \overline{\widehat{\psi}(\xi) |\xi|^{-\alpha}} d \xi \\
&=& \int_{R^n} \widehat{v_k}(\xi) \overline{\widehat{\psi}(\xi)} d \xi \\
&=& \int_{R^n} v_k(x) \psi (x) d x .
\end{eqnarray*}
Here we have used a result in \cite{LL} that the Fourier transform of $\int_{R^n} \frac{\psi(x)}{|x-y|^{n-\alpha}} d x$ is a constant multiple of $\widehat{\psi}(\xi) |\xi|^{-\alpha}$.
It follows that
\begin{eqnarray*}
\int_{R^n}g_k(x)\psi(x)dx &=& I_1 + I_2 \\
&=& \int_{R^n} \left( v_k(x) + \frac{c_1}{|x|^{n-\alpha}} \right) \psi (x) d x \\
&=& \int_{R^n} u_k(x) \psi (x) d x .
\end{eqnarray*}
This proves (\ref{18}) and hence (\ref{17}), and from which, we arrive immediately that

\begin{equation}
\int_{R^n\setminus B_r(0)}\frac{f_k(x)}{|x|^{n-\alpha}}dx \leq u_k(0) \leq u(0).
\label{20}
\end{equation}

Combining (\ref{7}) with (\ref{20}), we derive that
\begin{eqnarray}
I_2&=&c\int_{R^n\setminus
B_r(0)}\int_{R^n}\frac{u_k(x)-u_k(y)}{|x-y|^{n+\alpha}}dy \; \varphi(x)dx \nonumber\\
&=& c\int_{R^n\setminus
B_r(0)}f_k(x)\varphi(x)dx \nonumber\\
&\leq & c\int_{R^n\setminus
B_r(0)}\frac{f_k(x)}{|x|^{n-\alpha+1}}dx \nonumber\\
&\leq & \frac{c}{r}\int_{R^n\setminus
B_r(0)}\frac{f_k(x)}{|x|^{n-\alpha}}dx \nonumber\\
&\leq& \frac{c}{r} u(0) \rightarrow
0\,,\,\,\,as\,r\rightarrow \infty, \mbox{ uniformly in } k .  \label{21}
\end{eqnarray}

(\ref{15}) and (\ref{21}) imply that (\ref{11}) holds. Hence we have
\begin{eqnarray}
\lim_{k\rightarrow \infty}\int_{R^n}u_k(x)\psi(x)dx
&=&\lim_{k\rightarrow
\infty}\int_{R^n}u_k(x)(-\Delta)^{\alpha/2}\varphi(x)dx
\nonumber\\
&=&\lim_{k\rightarrow
\infty}\int_{R^n}(-\Delta)^{\alpha/2}u_k(x)\varphi(x)dx =0.
\label{22}
\end{eqnarray}

That is $$\int_{R^n}u(x)\psi(x)dx=0.$$

Therefore we come to the conclusion that $u\equiv C$.

This complete the proof of Theorem \ref{mthm2}.
\bigskip

{\bf The Proof of Theorem \ref{mthm1}.}

To see that Theorem \ref{mthm2} implies Theorem \ref{mthm1}, let $v$ be any $\alpha$-harmonic
function that is bounded from above by a constant $M$ in $R^n$. Take $u(x) = M - v(x)$, then
$$\left\{\begin{array}{ll}
(- \Delta)^{\alpha/2} u(x) = 0 , & x \in R^n , \\
u(x) \geq 0 , & x \in R^n .
\end{array}
\right.
$$
By Theorem \ref{mthm2}, $u$ must be constant, hence so does $v$. Similarly, if $v$ is any $\alpha$-harmonic function that is bounded from below by a constant $M$ in $R^n$, then we let $u(x) = v(x) - M$ to derive that $v$ must be constant. This completes the proof of Theorem \ref{mthm1}.

\section{Applications}

{\bf{The Proof of Theorem \ref{mthm3}.}}

Assume $u \in L_{\alpha}$ is a nonnegative locally bounded strong solution of
\begin{equation}
(- \Delta)^{\alpha/2}u(x)=u^p(x), \;\; x \in R^n,
\label{3.1}
\end{equation}

Let
\begin{equation}
v_R(x)=\int_{B_R}G_R(x,y)u^p(y)dy,\label{3-6}
\end{equation}
where $G_R(x,y)$ is Green's function on the ball $B_R(0)$:
\begin{eqnarray*}\left\{\begin{array}{ll}
(- \Delta)^{\alpha/2}G_R(x,y)=\delta(x-y), &  x, y \in B_R(0),\\
G_R(x,y)=0, & x \mbox { or } y \in  B^c_R(0).
\end{array} \right.\end{eqnarray*}
Thanks to Kulczycki \cite{Ku}, one can write
$$G_R(x,y)=\frac{A_{n,\alpha}}{s^{\frac{(n-\alpha)}{2}}}\left[1-\frac{B_{n,\alpha}}
{(s+t)^{\frac{(n-\alpha)}{2}}}\int^{\frac{s}{t}}_{0}\frac{(s-tb)^{\frac{(n-\alpha)}{2}}}{b^{\frac{\alpha}{2}}(1+b)}db \right],\,\,x,y\in B_R(0),$$
where $s=\frac{|x-y|^2}{R^2}$, $t=(1-\frac{|x|^2}{R^2})(1-\frac{|y|^2}{R^2})$. $A_{n,\alpha}$ and $B_{n,\alpha}$ are constants depending on $n$ and $\alpha$.

It is easy to verify that
\begin{eqnarray}\left\{\begin{array}{ll}
(- \Delta)^{\alpha/2}v_R(x)=u^p(x),& \mbox {in\,\,\,} B_R(0),\\
v_R=0,& \mbox { in\,\,\,} B^c_R(0).
\end{array} \right.\label{3-5} \end{eqnarray}

Let $w_R(x)=u(x)-v_R(x)$, by (\ref{3.1}) and (\ref{3-5}), we have
\begin{eqnarray}\left\{\begin{array}{ll}
(- \Delta)^{\alpha/2}w_R(x)=0,& \mbox {in\,\,\,} B_R(0),\\
w_R\geq 0,& \mbox { in\,\,\,} B^c_R(0).
\end{array} \right.\label{3-7} \end{eqnarray}

Applying the Maximum Principle \cite{Si} to (\ref{3-7}), we derive that
\begin{equation}
w_R(x)\geq 0,\,\,x\in R^n.\label{3-8}
\end{equation}

One can verify that,
\begin{equation}
v_R(x)\rightarrow v(x) = \int_{R^n}\frac{c_n}{|x-y|^{n-\alpha}}u^p(y)dy , \,\,\,as\,R\rightarrow\infty.\label{3-9}
\end{equation}

It's easy to see
\begin{equation}
(- \Delta)^{\alpha/2}v(x)=u^p(x),\,\,x\in R^n.\label{3-3}
\end{equation}

Denote
$$w(x)=u(x)-v(x).$$

Then by (\ref{3.1}), (\ref{3-3}), (\ref{3-8}), and (\ref{3-9}),  we have
\begin{eqnarray}\left\{\begin{array}{ll}
(- \Delta)^{\alpha/2}w(x)=0,& \mbox {in\,\,\,} R^n,\\
w\geq 0,& \mbox { in\,\,\,} R^n.
\end{array} \right.\label{3-10} \end{eqnarray}

From Theorem \ref{mthm2}, we derive that $w\equiv C$. Then obviously,
\begin{equation}
u(x)=w(x)+v(x)\geq C,\,\,x\in R^n.\label{3-11}
\end{equation}

Next, we show that $C=0$. Otherwise, if $C>0$, then
\begin{equation}
u(x) \geq v(x)=\int_{R^n}\frac{c_n}{|x-y|^{n-\alpha}}u^p(y)dy\geq \int_{R^n}\frac{c_nC^p}{|x-y|^{n-\alpha}}dy=\infty.\label{3-12}
\end{equation}

This is a contradiction.

Therefore we conclude that $$u(x)=v(x)=\int_{R^n}\frac{c_n}{|x-y|^{n-\alpha}}u^p(y)dy.$$

This complete the proof of Theorem \ref{mthm3}.
\bigskip

{\bf The Proof of Theorem \ref{mthm4}}.

It is a direct consequence of Theorem \ref{mthm3} and the following results from \cite{CLO} and \cite{CLO1}:

\begin{pro}
Assume that $n \geq 2$ and $u$ is a locally bounded nonnegative solution of the integral equation (\ref{3-12}) for $0<\alpha < n$. Then

i) In  the critical case when $p = \frac{n+\alpha}{n-\alpha}$, it must assume the form
$$ u(x) = c (\frac{t}{t^2 + |x -
x_o|^2})^{(n-\alpha)/2}$$
for some  $t>0$, $x_o\in{\mathbf{R}}^n$.

ii) In the subcritical case when $1<p<\frac{n+\alpha}{n-\alpha}$, we must have $u\equiv 0.$
\end{pro}

{\em Authors' Addresses and E-mails:}
\medskip

Ran Zhuo

Department of Mathematical Sciences

Yeshiva University

New York, NY, 10033 USA

ranzhuo201@gmail.com

\medskip

Wenxiong Chen

Department of Mathematical Sciences

Yeshiva University

New York, NY, 10033 USA

wchen@yu.edu

\medskip

Xuewei Cui

Department of Applied Mathematics

Northwestern Polytechnical University

Xian 710072, Shanxi, China

and Department of Mathematical Sciences

Yeshiva University

New York, NY, 10033 USA

xueweicui@hotmail.com

\medskip

Zixia Yuan

Department of Applied Mathematics

School of Mathematical Sciences

University of Electronic Science and Technology of China

No.2006, Xiyuan Ave, West Hi-Tech Zone, Chengdu, China

and Department of Mathematical Sciences

Yeshiva University

New York, NY, 10033 USA

yzx8047@yahoo.com.cn

\end{document}